\documentclass[12pt,reqno]{amsart}

\usepackage{amssymb}
\usepackage{latexsym}
\usepackage{amsmath}
\usepackage{amsthm}
\usepackage{amsxtra}

\usepackage{xcolor}
\usepackage{hyperref}

\usepackage{textgreek}

\usepackage{graphicx} 
\usepackage{float}
\usepackage{caption}
\usepackage{subcaption}

\title[Functional Iteration and Roots]{Analytical and Numerical Approaches\\ 
for Finding Functional Iterates and Roots}

\author[S. Nesargi]{Sanay Nesargi}
\address{Greenhill School, Addision, TX 75001}
\curraddr{}
\email{nesargis26@greenhill.org}

\author[G. Roudenko]{Gregory Roudenko}
\address{University of Florida,  Gainsville, FL 32612}
\curraddr{}
\email{groudenko@ufl.edu}

\date{}

\keywords{functional iteration, fractional iterates, functional roots, super-logarithm, tetration, numerical approximations}

\subjclass[2020]{39B05, 39B12, 39B22}

\begin{document}

\begin{abstract}
We investigate solutions to the functional equation \( f(f(x)) = e^{x} \), which can be interpreted as the problem of finding a half iterate of the exponential map. While no elementary solution exists, we construct and analyze non-elementary solutions using methods based on the Lambert $W$ function, tetration, and Abel’s functional equation. We examine structural properties of possible solutions, including monotonicity, injectivity, and behavior across different intervals, and provide a piecewise-defined framework that extends to the entire real domain.
Building on this, we introduce the super-logarithm and its inverse, the super-root, as analytic tools for defining fractional iterates of $e^{x}$. Using a power-series expansion near $x = 1$, we numerically approximate the super-logarithm and demonstrate a procedure for computing fractional iterates, including the half-iterate of the exponential function. Our approach is validated by comparisons to known constructions such as Kneser’s tetration, with an emphasis on computational feasibility and numerical stability.
Finally, we explore the broader landscape of fractional iteration, showing that similar techniques can be applied to other functions beyond $e^{x}$. Through numerical approximations and series-based methods using genetic algorithms and other optimization techniques, we confirm that fractional iterates not only exist for many analytic functions but can be computed with practical accuracy, opening pathways to further applications in dynamical systems and functional analysis.
\end{abstract}


\maketitle

\tableofcontents

\section{Introduction and Preliminaries}

The study of functional iteration concerns the repeated application of a function \(f\), resulting in compositions \(f^{\circ n}(x) = f(f(\cdots f(x)\cdots))\) applied \(n\) times. While the case of integer \(n\) is straightforward, the extension to fractional or continuous iterates — where \(n \in \mathbb{R}\) or even \(\mathbb{C}\) — poses significant challenges, e.g., see \cite{Erdos},  \cite{Chen}. Such generalizations require defining a function \(g\) such that \(g^{\circ n}\) interpolates smoothly between integer iterates while satisfying natural functional equations.

A particularly interesting instance of this problem arises in the search for a function \(f\) satisfying
\[
f(f(x)) = e^{x},
\]
which can be interpreted as finding a \emph{half iterate} or, more generally, the \emph{half-iterate} of the exponential map. In other words, we seek \(f\) such that applying it twice reproduces the usual exponential function. Despite the apparent simplicity of this equation, Section~2 shows that there is \emph{no elementary solution} to \(f(f(x)) = e^{x}\), using arguments from elementary function theory and the nonexistence of solutions expressible in standard closed-form operations.

To explore possible solutions, Section~3 analyzes the structural properties any function \(f\) must possess if it is to satisfy this functional equation, including \emph{injectivity}, \emph{monotonicity}, and constraints on the image of \(f\). We introduce a \emph{piecewise-defined solution}, extending it to different regions of the real line. Section~4 further considers related constructions under assumed forms, providing heuristic constructions that illustrate the general behavior of such functions.

The second half of this paper (Sections~5--10) shifts focus to tetration, the super-logarithm, and fractional iteration methods. Tetration, denoted \(x \uparrow\uparrow y\), generalizes repeated exponentiation to non-integer heights and provides the analytic framework needed for defining \(f^{\circ y}\) when \(y \in \mathbb{R}\). We introduce the super-logarithm \(\mathrm{slog}_x\), the super-root, and their relevance to Abel’s functional equation, which provides a general method for interpolating between iterates.

In Sections~6 and~7, we construct the \emph{numerical implementation} of the super-logarithm and tetration, beginning with a truncated power-series approximation near \(x = 1\) and extending it globally through functional iteration. We also present methods for numerical inversion and the computation of fractional iterates of the exponential function, including explicit approximations for the \emph{half-iterate}. Section~7.7 demonstrates this approach, verifying that fractional iterates do exist and can be computed to high accuracy, even though they cannot be expressed in elementary form.

Finally, Sections~8--10 explore broader \emph{numerical approximations}, including fixed points, Taylor series expansions, and periodic function iterations. We conclude by demonstrating that the methodology developed here applies not only to \(e^{x}\), but to a larger class of functions, showing that \emph{fractional iteration is both theoretically valid and numerically practical}. 

\subsection*{Acknowledgments}

Both authors would like to thank Stephen Wolfram for generously sponsoring the Wolfram Summer Research Program and Wolfram Emerging Leaders Program during which this paper was conceptualized. We are grateful Bertie Bennett and Christopher Gilbert for providing the original motivation for our research and the base Rust source code from which we developed our genetic algorithm. We also thank Ethan Dutton for providing consultation with generalized derivations using the Lambert $W$ function.

\section{Nonexistence of an Elementary Solution to 
$f(f(x)) = e^x$ }\label{nonexistence}

We begin by proving that no elementary function exists that satisfies the functional equation $f(f(x)) = e^x$.

Assume, for contradiction, that there exists an elementary function \( f(x) \) satisfying the functional equation
\[
f(f(x)) = e^x.
\]

Differentiating both sides with respect to \( x \), and applying the chain rule to the left-hand side, we obtain
\[
\frac{d}{dx} f(f(x)) = f'(f(x)) \cdot f'(x),
\]
and to the right-hand side,
\[
\frac{d}{dx} e^x = e^x.
\]
Hence,
\begin{equation}\label{E:1}
f'(f(x)) \cdot f'(x) = e^x.
\tag{1}
\end{equation}

\subsection{Nested iteration for \( f'(x) \)}

Rewriting equation \eqref{E:1}, we isolate \( f'(x) \):
\begin{equation}
f'(x) = \frac{e^x}{f'(f(x))}.
\tag{2}
\end{equation}

This defines $f'(x)$ in terms of $f'(f(x))$, which itself depends on $f'(f(f(x)))$, and so on. Recursively expanding this expression yields:
\[
f'(x) = \frac{e^x}{\dfrac{e^{f(x)}}{f'(f(f(x)))}} = \cdots,
\]
which becomes an infinitely nested, self-referential expression.

Such recursive derivative relations resemble earlier work on functional iteration methods for solving equations, where iterative descent often prevents closed-form elementary solutions from arising (see \cite{Hamilton1946}).

\subsection{Argument from Elementary Function Theory}

Elementary functions are defined as functions obtained from a finite number of compositions and combinations of rational functions, exponentials, logarithms, and trigonometric functions, along with their inverses. They are closed under:
\begin{itemize}
    \item Addition, subtraction, multiplication, division,
    \item Composition,
    \item Algebraic solutions of equations.
\end{itemize}

However, the recursive structure of equation (2) defines $f'(x)$ using an infinite descent through nested function evaluations. This process cannot terminate in a finite combination of elementary operations and functions. Therefore, the function $f(x)$ cannot be elementary.

\subsection{Simple Consequence} 

We conclude that no elementary function satisfies the equation
\[
f(f(x)) = e^x.
\]

\section{Preliminaries: Structural Properties of Continuous Solutions}

We now attempt to construct a non-elementary solution to $f(f(x))=e^x$, that is continuous and differentiable.

We begin by examining necessary properties of a continuous or a differentiable function \( f:\mathbb{R} \to \mathbb{R} \) satisfying
\[
f(f(x)) = e^x.
\]

\subsection{The image of \( f \) is not surjective}

To the contrary, suppose that $f(\mathbb{R}) = \mathbb{R}$. Then
$$
f(f(\mathbb{R})) = f(\mathbb{R}) = \mathbb{R},
$$
however, 
$$
f(f(\mathbb{R})) = e^{\mathbb{R}} = (0, \infty).
$$
This implies \( \mathbb{R} = (0, \infty) \), a contradiction. Therefore, $f$ is not surjective.

Since $f$ is continuous, its image must be an interval. Moreover, the equation $f(f(x)) = e^x$ implies
$$
\text{Im}(f(f)) = (0, \infty) \subseteq \text{Im}(f),
$$
so the image of $f$ must include all positive real numbers. It follows that
\[
\text{Im}(f) = (a, \infty) \quad \text{or} \quad [a, \infty) \quad \text{for some } a \leq 0.
\]
Thus, $f$ is bounded below.

\subsection{Injectivity of $f$}

Assume $f(y) = f(z)$. Then
\[
f(f(y)) = f(f(z)) \Rightarrow e^y = e^z \Rightarrow y = z,
\]
since the exponential function is injective. Therefore, $f$ itself is injective.
This aligns with general structural results on continuous solutions of functional equations involving iterates and powers, where injectivity often emerges as a necessary condition, e.g., see \cite{Morawiec2014}.

\subsection{Monotonicity of $f$}

Assuming differentiability, we differentiate both sides to get
\[
\frac{d}{dx} f(f(x)) = \frac{d}{dx} e^x.
\]
By the chain rule,
\[
f'(f(x)) \cdot f'(x) = e^x.
\]
Since \( e^x > 0 \), it follows that
\[
f'(x) \neq 0 \quad \text{for all } x.
\]
Hence, $f$ is strictly monotonic: either strictly increasing or strictly decreasing.
This conclusion is consistent with the general monotonicity restrictions derived in the study of iterative functional equations, \cite{Morawiec2014}.

\subsection{Summary}

Any differentiable solution $f$ to the equation $f(f(x)) = e^x$ satisfies the following:
\begin{itemize}
    \item $f$ is bounded below: $\exists a \leq 0$ such that $f(x) \geq a$ for all $x$,
    \item $f$ is injective,
    \item $f$ is strictly monotonic,
    \item $f$ is not surjective: \( f(\mathbb{R}) \subsetneq \mathbb{R} \). 
\end{itemize}

These conditions impose strong restrictions on any possible solution, but they do not uniquely determine \( f \). However, they help guide the construction of solutions.

\section{A Known Piecewise-Defined Solution}

Although no elementary closed-form solution exists for \( f(f(x)) = e^x \), a known solution can be constructed piecewise.

We seek a function \( f : \mathbb{R} \to \mathbb{R} \) satisfying the equation 
\[
f(f(x)) = e^x,
\]
defined piecewise using intervals where the behavior of the function is constructed recursively or approximated numerically.

Let us define \( f \) on an interval \( [a, b] \), then extend to \( \mathbb{R} \) by the functional equation. Specifically, we can define \( f \) on an initial interval \( I \), and then define it on the rest of the real line by setting
\[
f(x) = f^{-1}(e^x),
\]
or using iterated compositions and inverses depending on domain direction.

One known approach is to define \( f \) on an interval such as \( [0, 1] \) by some carefully chosen continuous and monotonic function \( \phi \), then extend it to \( \mathbb{R} \) recursively:
\[
f(x) =
\begin{cases}
\phi(x), & x \in [0, 1], \\
f^{-1}(e^x), & x > 1, \\
\log(f(f(x))), & x < 0.
\end{cases}
\]

This defines a function that satisfies \( f(f(x)) = e^x \) by construction, though it may not be smooth or differentiable everywhere unless further conditions are imposed on \( \phi \). The continuity and monotonicity of \( \phi \) can be adjusted to ensure global injectivity and boundedness.

We now show that the extensions of this function to \( x > 1 \) and \( x < 0 \) preserve the functional equation \( f(f(x)) = e^x \).

\subsection{Extend \( f \) to \( x > 1 \)}

For \( x > 1 \), define
\[
f(x) := f^{-1}(e^x).
\]
This definition is valid because \( f \) is strictly increasing on \([0,1]\), so the inverse \( f^{-1} \colon [1, e] \to [0,1] \) exists and is well-defined.

Furthermore, since \( x > 1 \Rightarrow e^x > e \), we recursively apply the same formula 
\[
f(x) = f^{-1}(e^x), \quad f(f(x)) = f(f^{-1}(e^x)) = e^x.
\]

\subsection{Extend \( f \) to \( x < 0 \)}

For \( x < 0 \), we define \( f \) recursively as
\[
f(x) := \log(f^{-1}(x)).
\]
We justify this by observing that
\[
f(f(x)) = f(\log(f^{-1}(x))) = e^{\log(f^{-1}(x))} = f^{-1}(x).
\]
Thus, applying \( f \) again, we obtain
\[
f(f(x)) = f^{-1}(x) \Rightarrow f(f(f(x))) = x \Rightarrow f(f(x)) = e^x \text{ if } x = \log(f(f(x))).
\]
Thus, defining \( f(x) = \log(f^{-1}(x)) \) in this way, ensures the same recursive property works backwards.

\subsection{Final Piecewise Definition}

All together:
\[
f(x) =
\begin{cases}
\phi(x) & \text{if } x \in [0,1], \\
f^{-1}(e^x) & \text{if } x > 1, \\
\log(f^{-1}(x)) & \text{if } x < 0.
\end{cases}
\]

\subsection{Verification}

We now verify that this function satisfies \( f(f(x)) = e^x \) for all \( x \in \mathbb{R} \):

\begin{itemize}
    \item For \( x \in [0,1] \), we have \( f(f(x)) = f(\phi(x)) = e^x \) by construction.
    \item For \( x > 1 \), \( f(x) = f^{-1}(e^x) \Rightarrow f(f(x)) = f(f^{-1}(e^x)) = e^x \).
    \item For \( x < 0 \), \( f(x) = \log(f^{-1}(x)) \Rightarrow f(f(x)) = f(\log(f^{-1}(x))) = e^{\log(f^{-1}(x))} = f^{-1}(x) \). Then \( f(f(x)) = e^x \) by the inverse property.
\end{itemize}

\subsection{Final Definition}

This piecewise construction defines a function \( f \) that is:
\begin{itemize}
    \item Continuous (assuming a smooth choice of \( \phi \) on \([0,1]\)),
    \item Strictly increasing on its domain,
    \item Satisfying \( f(f(x)) = e^x \) for all \( x \in \mathbb{R} \).
\end{itemize}

\section{Solutions to a Related Functional Equation Under an Assumed Form} \label{assumedForm}

We now start by setting the following equal to \( e^x \):
\[
y^{y^x} = e^x.
\]

Taking the natural logarithm of both sides, it is straightforward to get
\[
\log(y^{y^x}) = \log(e^x) \Rightarrow y^x \cdot \log y = x.
\]

Let \( \log y = u \), so \( y = e^u \). Then
\[
y^x = e^{ux} \quad \text{and} \quad y^x \log y = e^{ux} \cdot u = x.
\]

This gives us that
\[
u e^{ux} = x.
\]

Let \( z = ux \), so \( u = \frac{z}{x} \). Substituting, we get
\[
\frac{z}{x} e^z = x \Rightarrow z e^z = x^2.
\]

By definition of the Lambert $W$ function, \( z = W(x^2) \). Therefore,
\[
u = \frac{W(x^2)}{x} \Rightarrow y = e^u = e^{W(x^2)/x}.
\]

This is equivalent to \( \left( \frac{x^2}{W(x^2)} \right)^{1/x} \), since
\[
\left( \frac{x^2}{W(x^2)} \right)^\frac{1}{x} = \left(e^{W(x^2)}\right)^{\frac{1}{x}}.
\]

Again, by definition of the Lambert $W$ function, \( W(x^2)e^{W(x^2)} = x^2 \), therefore, 
\[
e^{W(x^2)} = \frac{x^2}{W(x^2)} \Rightarrow y = \left( \frac{x^2}{W(x^2)} \right)^{1/x}
\]

Therefore, the function
\[
\boxed{f(x) = \left( \frac{x^2}{W(x^2)} \right)^{1/x}}
\]
satisfies \( f(x)^{f(x)^{x}} = e^x \) for all \( x > 0 \).

While this function is not elementary, it is analytic and well-defined for \( x > 0 \). The result demonstrates that even without a general closed-form approach, making a strategic assumption about the structure of \( f \) can yield a valid solution—at least over a restricted domain. This opens the door to further exploration of uniqueness, smooth extensions, or possible generalizations beyond the positive real axis.

\subsection{Generalization}
This same technique used above in this section \ref{assumedForm} can be used for similar functions in the exponential family. A prime example is $xe^x$, where by following the same steps as earlier one can arrive at 
$$
\bigg( \frac{x\ln{x}+x^2}{W(x\ln{x}+x^2)} \bigg)^\frac{1}{x}.
$$

It is likewise possible to extend the above method to all functions of the form $E(x) = a(x)\cdot{}e^{b(x)}$. 
The form of the function $y$ then becomes
$$
\bigg(\frac{x\ln{a(x)}+x^2b(x)}{W(x\ln{a(x)}+x^2b(x))} \bigg)^\frac{1}{x}.
$$

Therefore, we have 
\[\boxed{f(x)=\frac{x\ln{a(x)}+x^2b(x)}{W(x\ln{a(x)}+x^2b(x))}}\]

\section{Tetration, the Super-Logarithm, and the Super-Root}

In the hierarchy of arithmetic operations, each operation can be seen as repeated application of the previous:
\begin{itemize}
    \item Addition: \( a + b \)
    \item Multiplication: repeated addition, \( a \times b = \underbrace{a + a + \dots + a}_{b \text{ times}} \)
    \item Exponentiation: repeated multiplication, \( a^b = \underbrace{a \times a \times \dots \times a}_{b \text{ times}} \)
    \item \textbf{Tetration}: repeated exponentiation, \( ^b a = \underbrace{a^{a^{\cdot^{\cdot^{a}}}}}_{b \text{ times}}. \)
\end{itemize}

Formally, tetration is defined recursively for natural \( n \) by:
\[
^1 a = a, \quad {}^{n+1}a = a^{({}^n a)} \quad \text{for } n \geq 1.
\]
It grows at an extremely fast rate. For example:
\[
{}^1 2 = 2, \quad {}^2 2 = 4, \quad {}^3 2 = 16, \quad {}^4 2 = 65536.
\]

\subsection{The Super-Logarithm}

Just as the logarithm is the inverse of exponentiation, the \textbf{super-logarithm} (denoted \( \operatorname{slog}_a(x) \)) is the inverse of tetration:
\[
\operatorname{slog}_a(x) = y \quad \Longleftrightarrow \quad {}^y a = x
\]
This function can be thought of as ``how many times must we exponentiate \( a \) to get \( x \)?"

\subsection{The Super-Root}

The \textbf{super-root} is the inverse operation in the base rather than the height. It is defined as:
\[
\operatorname{sroot}_n(x) = a \quad \Longleftrightarrow \quad {}^n a = x
\]
In other words, it answers the question: "what base \( a \) tetrated \( n \) times gives \( x \)?"

\subsection{Relevance to the Functional Equation}

The equation \( f(f(x)) = e^x \) can be interpreted as a \emph{half iterate} of the exponential:
\[
f \circ f = \exp.
\]
This is conceptually similar to asking: “What function, when applied twice, gives exponentiation?” — analogous to asking what number squared gives another (i.e., square root), but lifted into the function space.

This invites the use of operations like tetration and its inverses. For instance, define `$slog$'
\[
f(x) = \operatorname{slog}_a(e^x)
\]
Then
\[
f(f(x)) = \operatorname{slog}_a(e^{\operatorname{slog}_a(e^x)}) = \operatorname{slog}_a({}^a(e^x)) = x
\]
(assuming suitable domains and definitions of \(\operatorname{slog}_a\) and its inverse). Alternatively, we may define
\[
f(x) = \operatorname{sroot}_2(e^x),
\]
if such a function satisfies
\[
f(f(x)) = e^x \quad \Longleftrightarrow \quad {}^2 f(x) = e^x.
\]

Hence, tetration-based functions provide a conceptual framework for expressing or approximating solutions to this type of recursive functional equation, even if closed forms are elusive. They also open the door to extending the notion of function iteration beyond integers and exploring fractional or even continuous function iteration — a major area of research in functional equations and dynamical systems.

\section{Hyperoperation-based Solutions}

\subsection{Tetration and Its Recursive Property}

Tetration is the next hyperoperation after exponentiation and is denoted as
\[
x \uparrow\uparrow y : = \underbrace{x^{x^{\cdot^{\cdot^{x}}}}}_{y \text{ times}}.
\]
It satisfies the fundamental recursive property:
\[
x \uparrow\uparrow y = x^{(x \uparrow\uparrow (y-1))},
\]
with base case:
\[
x \uparrow\uparrow 1 = x.
\]
For example, 
\[
x \uparrow\uparrow 2 = x^x, \quad x \uparrow\uparrow 3 = x^{x^x}.
\]

\subsection{The Super-Logarithm}

The super-logarithm, denoted as \( \mathrm{slog}_x \), is defined as the inverse of tetration:
\[
\mathrm{slog}_x(x \uparrow\uparrow y) = y,
\]
which implies:
\[
x \uparrow\uparrow \mathrm{slog}_x(z) = z.
\]

Using the recursive property of tetration, we can derive a similar recursive property for the super-logarithm.  
Starting from
\[
x \uparrow\uparrow y = x^{(x \uparrow\uparrow (y-1))},
\]
taking \( \mathrm{slog}_x \) of both sides gives
\[
\mathrm{slog}_x(x \uparrow\uparrow y) = \mathrm{slog}_x(x^{(x \uparrow\uparrow (y-1))}) = \mathrm{slog}_x(x \uparrow\uparrow (y-1)) + 1.
\]
This generalizes to the identity
\[
\mathrm{slog}_x(x^z) = \mathrm{slog}_x(z) + 1,
\]
which we refer to as Property (2).

\subsection{Connection to Abel’s Functional Equation}

A key observation is that Property (2) mirrors Abel's functional equation:
\[
\psi(f(x)) = \psi(x) + 1,
\]
where \( \psi \) is a linearizing function for \( f \). From Abel’s equation, it follows that
\[
\psi(f^n(x)) = \psi(x) + n,
\]
and hence,
\[
f^{n}(x) = \psi^{-1}(\psi(x) + n),
\]
for any real \( n \), allowing us to define continuous (even fractional) iterates of \( f \).

For exponentials, we set
\[
f(z) = x^z, \quad \psi(z) = \mathrm{slog}_x(z),
\]
which immediately gives
\[
\mathrm{slog}_x(x^z) = \mathrm{slog}_x(z) + 1,
\]
matching Property (2). Therefore, tetration and the super-logarithm naturally emerge as inverse functions that linearize the iteration of exponentials.

\subsection{Analytic Definition of the Super-Logarithm}

To be analytically useful, we require \( \mathrm{slog}_x \) to be infinitely differentiable (i.e., \( C^{\infty} \)) and analytic. Such constructions are detailed in \cite{tetration}, where \( \mathrm{slog}_x \) is defined to satisfy
\[
\mathrm{slog}_x(x \uparrow\uparrow y) = y, \quad \mathrm{slog}_x(x^z) = \mathrm{slog}_x(z) + 1,
\]
with smooth interpolation between integer values of \( y \).

\subsection{Fractional Iterates of the Exponential Function}

Using Abel’s functional equation and the super-logarithm, we can define the \( n \)-th iterate of \( \exp \) (even for non-integer \( n \)) as:
\[
\exp^{[n]}(x) = \mathrm{tet}_e \big( \mathrm{slog}_e(x) + n \big),
\]
where \( \mathrm{tet}_e \) is the inverse of \( \mathrm{slog}_e \), i.e., tetration to base \( e \).

\subsection{The Half-Iterate of the Exponential}

For \( n = \frac{1}{2} \), we define:
\[
\exp^{[1/2]}(x) = \mathrm{tet}_e \big( \mathrm{slog}_e(x) + \frac{1}{2} \big),
\]
which is the unique analytic function that satisfies:
\[
\exp^{[1/2]}(\exp^{[1/2]}(x)) = \exp(x).
\]
While this function has no elementary closed form, it can be approximated numerically. For example:
\[
\exp^{[1/2]}(1) \approx 1.291285997.
\]
This value can be obtained by numerically solving for the function \( f \) that satisfies:
\[
f(f(1)) = e,
\]
where \( f \) is smooth and monotonic, and then evaluating \( f(1) \).

\section{Numerical Implementation of the Super-Logarithm and Tetration}
\label{sec:implementation}

\subsection{Overview}
To numerically approximate real-height tetration $\operatorname{tet}_b(y)$ and its inverse, the super-logarithm $\operatorname{slog}_b(z)$, we adopt a hybrid local--global approach. A truncated power series representation of $\operatorname{slog}_b$ is constructed near $z=1$, with coefficients determined by enforcing the defining functional equation
\begin{equation}
\operatorname{slog}_b(b^{z}) = \operatorname{slog}_b(z) + 1.
\label{eq:slog-def}
\end{equation}
This local model is then extended globally by iterative application of logarithms or exponentials, while numerical inversion of the local series via root-finding allows the evaluation of $\operatorname{tet}_b(y)$ for non-integer heights.

\subsection{Series Construction and the Coefficient Matrix}
Locally, the super-logarithm is approximated by the truncated series
\begin{equation}
S(z) = -1 + \sum_{k=1}^{n} \frac{c_k}{k!}\, z^{k},
\label{eq:slog-series}
\end{equation}
where $n$ is the truncation order and $c_k$ are unknown coefficients. To determine the vector 
$\mathbf{c} = (c_1, \ldots, c_n)$, we linearize \eqref{eq:slog-def} around $z=1$ and enforce the equality of the first $n$ derivatives on both sides.

This leads to a system of $n$ linear equations,
\begin{equation}
\sum_{k=1}^{n} \left( \frac{k^{j}}{k!} - \delta_{j,k} (\log b)^{-k} \right) c_k =
\begin{cases}
1, & j = 0, \\
0, & 1 \le j \le n-1,
\end{cases}
\label{eq:linear-system}
\end{equation}
where $\delta_{j,k}$ is the Kronecker delta,
\begin{equation}
\delta_{j,k} =
\begin{cases}
1 & \text{if } j=k, \\
0 & \text{otherwise}.
\end{cases}
\end{equation}
The Kronecker delta ensures that the term $(\log b)^{-k}$ is subtracted only along the diagonal of the coefficient matrix. In matrix form, \eqref{eq:linear-system} reads
\begin{equation}
M \mathbf{c} = \mathbf{e}_1, \quad
M_{j+1,k} = \frac{k^{j}}{k!} - \delta_{j,k} (\log b)^{-k},
\end{equation}
with $\mathbf{e}_1 = (1,0,\ldots,0)^{T}$. The Mathematica routine
\texttt{SuperLogPrepare[n, b]} constructs $M$ and solves for $\mathbf{c}$ using \texttt{LinearSolve}.

\subsection{Global Evaluation of $\operatorname{slog}_b(z)$}

The polynomial approximation $S(z)$ is accurate only for $z \in (0,1]$. For arguments $z > 1$, we repeatedly apply the base-$b$ logarithm,

\begin{equation}
z \mapsto \log_b z = \frac{\ln z}{\ln b},
\end{equation}
until the iterated value lies within $(0,1]$. If $i$ iterations are applied, the global super-logarithm is approximated by
\begin{equation}
\operatorname{slog}_b(z) \approx S(\log_b^{\circ i}(z)) + i.
\end{equation}
For $z \le 0$, a symmetric procedure using exponentiation,
\[
z \mapsto b^z,
\]
is applied to shift the argument into $(0,1]$, with a correction term of $-1$. This process is implemented in the function \texttt{SuperLog}.

\subsection{Numerical Inversion and Tetration}
To evaluate $\operatorname{tet}_b(y)$ for real $y$, we decompose $y$ into an integer and a fractional part,
\[
y = m + r, \quad m = \lceil y \rceil, \quad r = y - m \in (-1,0].
\]
We first solve for $z$ in
\[
S(z) = r,
\]
which is accomplished numerically using Newton's method via \texttt{FindRoot}, initialized near $z=1$. We then reconstruct the integer component by iterating the base map:
\begin{equation}
\operatorname{tet}_b(y) \approx
\begin{cases}
b^{b^{\cdot^{\cdot^{b^{z}}}}}, & \text{($m$ times), } y > -1, \\
\log_b^{\circ (-m)} z, & y \le -1,
\end{cases}
\end{equation}
where $z = \text{TetCrit}(r)$. This procedure is implemented by the function \texttt{Tetrate}.

\subsection{Relation to Kneser's Solution}
Our implementation provides a practical and numerically efficient approximation of $\operatorname{slog}_b$ on the real line. However, it is fundamentally local and non-analytic due to the truncation in \eqref{eq:slog-series}. In contrast, Kneser \cite{kneser1950} constructed a globally analytic continuation of tetration for real bases $b > e^{1/e}$ by solving the Abel functional equation
\[
\phi(b^{x}) = \phi(x) + 1,
\]
where $\phi(x)$ (the Schröder function) \cite{Schroder} is analytically continued into the complex plane. The tetration function is then defined as
\[
F(z) = \phi^{-1}(\phi(1) + z),
\]
which satisfies $F(z+1) = b^{F(z)}$ exactly and is holomorphic on $\mathbb{C} \setminus (-\infty, 0]$. While Kneser's method provides the canonical analytic tetration, it is computationally expensive, requiring complex contour integration and high-precision evaluation. Our series-based approach can thus be viewed as a low-cost, real-axis approximation of Kneser's function, sufficient for experimental and numerical purposes.

\subsection{Numerical Stability}
The accuracy of this approach depends on both the truncation order $n$ and the conditioning of the matrix $M$. Larger $n$ improves local accuracy of $S(z)$ but can introduce ill-conditioning, requiring arbitrary-precision arithmetic. The Newton-based inversion (\texttt{FindRoot}) is sensitive to the quality of the initial guess and the validity range of the series $S(z)$.

\subsection{Fractional Iteration of the Exponential Function}

The primary goal of our construction is to enable not only integer iterates of the exponential map
\begin{equation}
E_b(x) = b^{x},
\end{equation}
but also fractional iterates $E_b^{\circ n}(x)$, where $n \in \mathbb{R}$ (or even $n \in \mathbb{C}$). Abel's functional equation \cite{Sz} provides a theoretical framework for this extension:
\begin{equation}
\Phi(E_b(x)) = \Phi(x) + 1,
\label{eq:abel}
\end{equation}
where $\Phi(x)$ is an invertible function (analogous to our $\operatorname{slog}_b(x)$). The $n$-th iterate of $E_b$ can then be defined by
\begin{equation}
E_b^{\circ n}(x) = \Phi^{-1}(\Phi(x) + n).
\label{eq:abel-iterate}
\end{equation}

Our numerical implementation follows this principle directly. Given the environment vector \texttt{env} computed by \texttt{SuperLogPrepare}, we define:
\begin{verbatim}
IntOrFractionalIterate[n_, x_] := Tetrate[env, SuperLog[env, x] + n]
\end{verbatim}
This function first maps $x$ into the $\Phi$-coordinate system using $\operatorname{slog}_b(x) = \texttt{SuperLog[env, x]}$, then shifts by $n$, and finally maps back to the original domain using $\operatorname{tet}_b$.

For instance, the case $n = \tfrac{1}{2}$ corresponds to the "half iterate" of the exponential function, formally defined as
\[
E_b^{\circ \frac{1}{2}}(x) = \operatorname{tet}_b \!\left( \operatorname{slog}_b(x) + \frac{1}{2} \right),
\]
which interpolates smoothly between $x$ and $E_b(x)$.

In code, we define:
\begin{verbatim}
nIt = 1/2;
IOFI[x_] := Tetrate[env, SuperLog[env, x] + nIt]
\end{verbatim}
Thus, \texttt{IOFI[x]} computes the fractional (half) iterate $E_b^{\circ 1/2}(x)$. By varying \texttt{nIt}, one can compute arbitrary real (or fractional) iterates of the exponential function, effectively solving the original functional iteration problem.

\section{Fixed points and Taylor Series}

It is possible to construct a Taylor series to approximate a half iterate by assuming a fixed point.

sin(x) has a fixed point at (0,0). Suppose the half iterate of sin, "$rin(x)$" also has a fixed point here. We can then construct the following functional equations by differentiating the definition of $rin(x)$.

\[sin(x)=rin(rin(x))\]
\[cos(x) = rin'(rin(x))rin'(x)\]
\[-sin(x) = rin'(x)^2rin''(rin(x))rin''(x)\]
\[ ... \]

Plugging in 0 for $x$ gives us

\[sin(0)=rin(rin(0)) \rightarrow rin(0) = 0\]
\[cos(0) = rin'(rin(0))rin'(0) \rightarrow rin'(0) = 1\]
\[-sin(0) = rin'(0)^2rin''(rin(0))rin''(0) \rightarrow rin''(0) = 0\]
\[ ... \]

Hence, using the value of each sequential derivative of $rin(x)$ at 0 we are able to construct a MacLaurin series for $rin(x)$.

This was done in \cite{maclaurinseries} and \cite{Chen}, with the result shown here. One can see that the approximation obtained using this result begins to diverge from the expected function outside the range $(-\pi/2,\pi/2)$, as this is the radius of convergence of the power series.

\begin{figure}[h]
    \centering
    \includegraphics[height=0.4\linewidth]{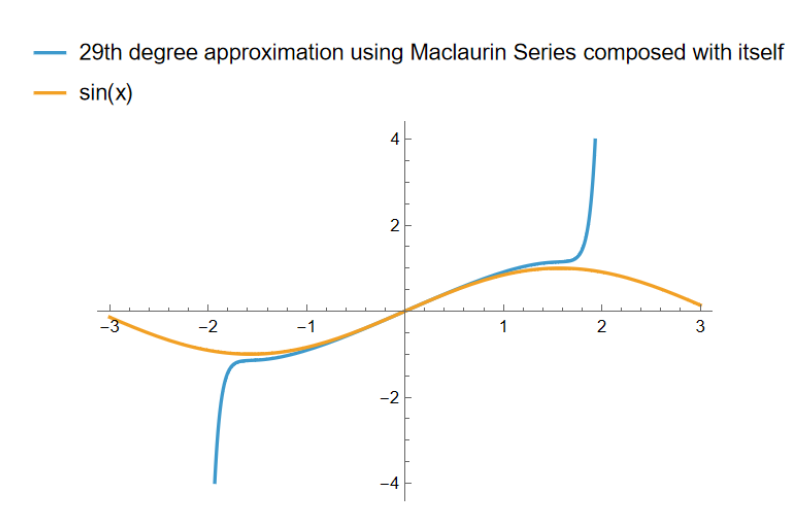}
    \caption{Approximation of $sin(x)$ using a Maclaurin series constructed with the fixed point method}
    \label{maclaurinseries}
\end{figure}

\section{Numerical Approximations of Functional  Roots}

\subsection{Introduction}

Although analytically finding iterates can prove useful and an interesting exercise, it is often difficult. In fact, in some cases it can be proven that it is impossible .

The problem of finding a half iterate $g$ of a function $f$,
$$
g(g(x)) = f(x),
$$
is a classical challenge in functional iteration. Except for certain simple cases, closed-form solutions for $g$ rarely exist, \cite{Hamilton1946}. However, in practical settings, we often do not need an exact expression for $g$ — a sufficiently accurate numerical approximation is enough. Modern approaches emphasize such approximation strategies for transcendental equations through iterative functional root algorithms \cite{Pakdemirli2024}.

The goal is to approximate $g$ on a chosen interval $[a,b]$ such that the error
$$
E(g) = \int_a^b \big(g(g(x)) - f(x)\big)^2 \, dx
$$
is minimized. To do so, we consider $g$ as an unknown function that we can approximate using flexible representations (e.g., piecewise-linear splines, polynomials, or Fourier series). We then fit the approximation by minimizing $E(g)$.

The following sections describe several numerical methods for approximating half iterates, starting with the Iterative Chain Approximation.

\subsection{Iterative Chain Approximation (ICA)}

\paragraph{\underline{Core Idea}}
The key observation is that if $g$ satisfies $g(g(x)) = f(x)$, then repeated application of $f$ generates constraints on $g$. Starting from an initial point $x_i$ and an initial guess $w$ for $g(x_i)$, we can define a sequence of points that $g$ must interpolate:
$$
(x_i, w), \, (w, f(x_i)), \, (f(x_i), f(w)), \, (f(w), f(f(x_i))), \, \ldots
$$
These pairs $(x_k, y_k)$ form a \emph{chain of constraints} that allow us to approximate $g$.
\smallskip

\paragraph{\underline{Construction of $g$}}
We proceed as follows:
\begin{enumerate}
    \item Generate sequences
    $$
    x_0 = x_i, \quad x_{k+1} = f(x_k),
    $$
    and
    $$
    y_0 = w, \quad y_{k+1} = f(y_k),
    $$
    for $k = 0,1,\ldots,n-1$.
    \item Use the set of points $\{(x_k, y_k)\}_{k=0}^n$ as interpolation data.
    \item Construct $g$ as a piecewise-linear spline through these points.
\end{enumerate}

\paragraph{\underline{Error Metric}}
To measure the quality of the approximation, we define
$$
E(g) \approx \sum_{x \in \mathcal{G}} \big(g(g(x)) - f(x)\big)^2 \, \Delta x,
$$
where $\mathcal{G}$ is a uniform grid in $[a,b]$.
\smallskip

\paragraph{\underline{Grid Search Optimization}}
The choice of $(x_i, w)$ strongly affects $g$. We therefore:
\begin{enumerate}
    \item Perform a coarse grid search over $(x_i, w)$ within a specified range.
    \item For each pair, compute $E(g)$ and record the best $(x_i, w)$.
    \item Refine the search by narrowing the range and reducing the step size.
    \item Repeat until convergence.
\end{enumerate}
\medskip

\paragraph{\underline{Pseudocode}}

\begin{verbatim}
Input: Function f(x), initial ranges for x_i and w

For refinement step = 1 to N:
    For x_i in [x_min, x_max] with step Δx:
        For w in [w_min, w_max] with step Δw:
            1. Generate chain (x_k, y_k) via f
            2. Build spline g(x) interpolating (x_k, y_k)
            3. Compute error
            4. If E < best_error: update best parameters
    Shrink [x_min, x_max] and [w_min, w_max] around the best pair
    Halve Δx and Δw

Return best g(x), x_i*, w*, and final error
\end{verbatim}

\paragraph{\underline{Implementation}}
A Python implementation of ICA (see \texttt{ICA.py}) generates chains, constructs linear splines using \texttt{scipy.interpolate}, and plots the approximation using \texttt{matplotlib}. A detailed description of the steps is given in the associated \texttt{ICA.md} documentation.

\subsection{Functional Root Refinement by Additive Correction}

\paragraph{\underline{Problem Setup}}
We continue solving the functional root equation:
$$
g(g(x)) = f(x),
$$
where $f$ is known, and $g$ (the "half iterate") is to be approximated.
While the Iterative Chain Approximation (ICA) provides an initial guess for $g$, this second method directly \emph{refines} an existing approximation by correcting its functional error.
\smallskip

\paragraph{\underline{Core Idea}}
Define the functional error:
$$
\delta(x) = f(x) - g(g(x)).
$$
If $\delta(x)$ vanishes for all $x$, then $g$ is an exact solution. Otherwise, we can \emph{update} $g$ by adding a scaled correction:
$$
g_{\text{new}}(x) = g(x) + \tau \, \delta(x),
$$
where $\tau$ (the "temperature") is a small positive constant that controls the adjustment step. 
This update moves $g$ in the direction that reduces the mismatch $g(g(x)) - f(x)$, similar to a gradient-free relaxation step.
\smallskip

\paragraph{\underline{Numerical Representation}}
To work numerically, we represent $g$ as a set of sample values:
\[
x_{\text{vals}} = \{x_0, x_1, \ldots, x_N\},
\]
\[
f_{\text{vals}} = f(x_{\text{vals}}),
\]
\[
g_{\text{vals}} = g(x_{\text{vals}}).
\]
A linear spline is built from $(x_{\text{vals}}, g_{\text{vals}})$ to define $g$ over the entire interval.

\paragraph{\underline{Algorithm}}
Starting from an initial guess $g_0$ (such as ICA’s result), we iterate:
\begin{enumerate}
    \item Evaluate $g(g(x))$ on the grid: $a_{\text{vals}} = g(g(x_{\text{vals}}))$.
    \item Compute $\delta = f_{\text{vals}} - a_{\text{vals}}$.
    \item Update $g$: $g_{\text{vals}} \leftarrow g_{\text{vals}} + \tau \delta$.
    \item Rebuild $g$ as a spline through the updated points.
    \item Evaluate the error
    $$
    E(g) = \int_{x_{\min}}^{x_{\max}} (g(g(x)) - f(x))^2 \, dx,
    $$
    and store $g$ if $E(g)$ improves.
\end{enumerate}

\paragraph{\underline{Pseudocode}}
\begin{verbatim}
Input: Function f(x), initial approximation g0(x)
x_vals = grid over [x_min, x_max]
f_vals = f(x_vals)
g = g0
for iteration = 1 to N:
    approx_vals = g(g(x_vals))
    delta = f_vals - approx_vals
    g_vals = g_vals + τ * delta
    rebuild spline g(x) from g_vals
    compute E(g)
return g with smallest E(g)
\end{verbatim}

\paragraph{\underline{Implementation}}
A Python implementation of Additive Correction (see \texttt{AdditiveCorrection.py}) performs the above algorithm, and plots the approximation using \texttt{matplotlib}. A detailed description of the steps is given in the associated \texttt{AdditiveCorrection.md} documentation.

\subsection{Genetic Algorithm of Fourier and Taylor Series Approximations}

\subsubsection{Overview}
This method approximates the \textit{half iterate} $g$ of a function $f$, where
\[
g(g(x)) = f(x).
\]
Instead of directly constructing $g$ through pointwise constraints, we assume that $g$ belongs to a \textit{parametric family of functions}, such as
\begin{itemize}
    \item \textbf{Fourier Series:} $\displaystyle g(x) = \sum_{k=1}^{n} a_k \sin(kx)$,
    \item \textbf{Taylor Series:} $\displaystyle g(x) = \sum_{k=0}^{n} \frac{c_k x^k}{k!}$,
    \item \textbf{Linear Splines:} $g$ is defined by linearly interpolating between control points.
\end{itemize}
We then \textit{optimize the coefficients} of the chosen family using a \textit{genetic algorithm} that evolves candidate solutions over multiple generations.

\subsubsection{Loss Function}
To evaluate how well a candidate $g$ satisfies the functional equation, we use the error defined in \label{12}.
The objective of the genetic algorithm is to minimize $\mathcal{L}(g)$.

\subsubsection{Genetic Algorithm Steps}
The algorithm proceeds as follows:
\begin{enumerate}
    \item \textbf{Initialization:} Generate a population of $N$ random candidate functions by randomly sampling their coefficients.
    \item \textbf{Evaluation:} Compute the loss $\mathcal{L}(g)$ for each candidate.
    \item \textbf{Selection:} Retain the top-performing candidates (lowest $\mathcal{L}$).
    \item \textbf{Mutation:} Create new candidates by adding random perturbations to the coefficients of the selected candidates. The mutation strength (``temperature'') decreases over time.
    \item \textbf{Iteration:} Repeat evaluation, selection, and mutation over many generations, improving the approximation of $g$.
\end{enumerate}

\subsubsection{Fourier and Taylor Series}
\paragraph{Fourier Series.}  
For periodic functions like $\sin x$, $g$ can be represented as
\[
g(x) = \sum_{k=1}^{n} a_k \sin(kx),
\]
where a small number of terms often yields a good approximation due to the series’ ability to capture oscillatory behavior.

\paragraph{Taylor Series.}  
For smooth, non-periodic functions, $g$ can be expressed as
\[
g(x) = \sum_{k=0}^{n} \frac{c_k x^k}{k!},
\]
which works well for approximations on finite intervals.

\subsubsection{Advantages and Limitations}
\begin{itemize}
    \item \textbf{Advantages:}
        \begin{itemize}
            \item Flexible approach, works for periodic and non-periodic functions.
            \item Capable of achieving very low error (e.g., for $\sin x$, $\mathcal{L}(g) < 10^{-12}$).
            \item Can switch between Fourier, Taylor, or spline-based representations.
        \end{itemize}
    \item \textbf{Limitations:}
        \begin{itemize}
            \item Slow convergence for high-degree series due to the large search space.
            \item Can get stuck in local minima, requiring multiple runs.
            \item Sensitive to tuning parameters such as mutation strength, population size, and series length.
        \end{itemize}
\end{itemize}

\subsubsection{Example Results}
\begin{itemize}
    \item A 9-term Fourier series for $f(x) = \sin x$ achieved a loss of $< 10^{-7}$.
    \item For $f(x) = 1 + x^2$, a Taylor polynomial yielded $\mathcal{L}(g) \approx 4.08 \times 10^{-5}$.
    \item For complex or piecewise functions, linear splines often outperform polynomials.
\end{itemize}

\subsubsection{Implementation}
A Python implementation of the genetic algorithm is based on our approach in \cite{WELPIterates} (see \texttt{genetic\_root\_solver.py}). This implementation represents $g$ as a parametric function (Fourier series, Taylor polynomial, or linear spline) and evolves its coefficients through mutation and selection. The loss $\mathcal{L}(g)$ is approximated using a discrete Riemann sum over the interval of interest.

The script leverages \texttt{numpy} for vectorized evaluation of candidate functions, \texttt{dataclasses} for managing coefficient sets, and \texttt{matplotlib} for visualizing $f(x)$ versus $g(g(x))$. A detailed description of the steps is provided in the associated documentation.

\begin{figure}[h]
    \centering
    \includegraphics[height=0.4\linewidth]{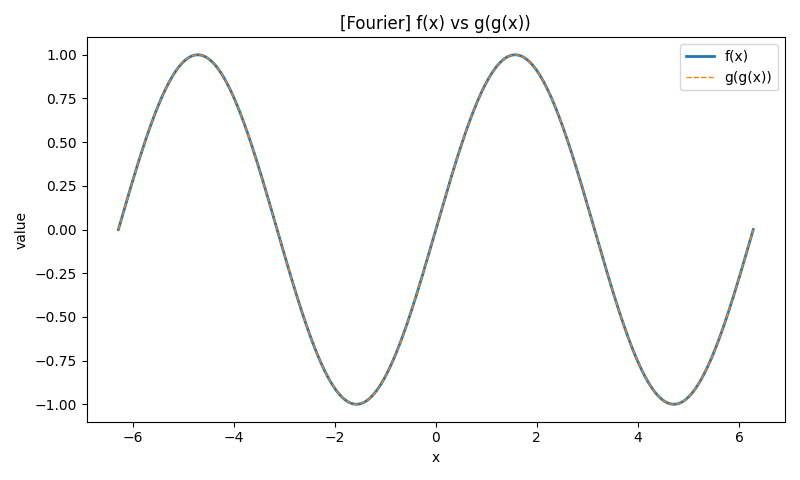}
    \caption{Genetic Algorithm with Fourier Series, approximating $\sin(x)$.}
    \label{geneticfourier}

    \centering
    \includegraphics[height=0.4\linewidth]{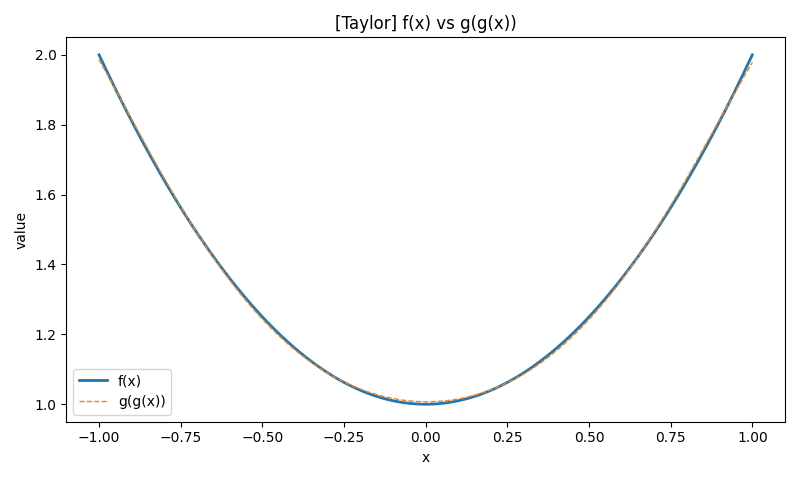}
    \caption{Genetic Algorithm with Taylor Series, approximating $x^2+1$.}
    \label{genetictaylor}
\end{figure}

\begin{figure}[h]
    \centering
    \includegraphics[height=0.4\linewidth]{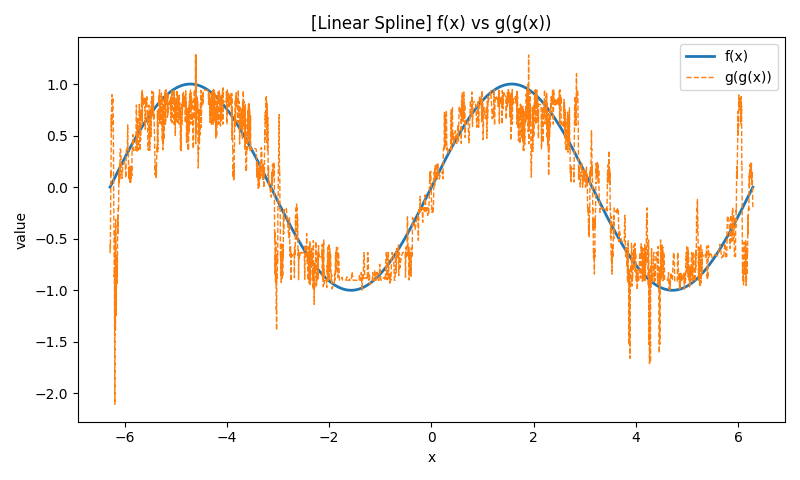}
    \caption{Genetic Algorithm using points of a linear spline, approximating $\sin(x)$.}
    \label{geneticlinspline}
\end{figure}

Shown in figures \ref{geneticfourier}, \ref{genetictaylor}, and \ref{geneticlinspline} are examples of approximations derived using this technique. Notice that using series of functional families to approximate gives a much more accurate and clean result than directly modifying points of a spline.

\section{Higher Order Approximations}

Using the genetic algorithm technique referenced in the previous section, it is possible to approximate any order iterate of a continuous function that resides within a certain family of functions and can be approximated as a linear combination of them (e.g. fourier/taylor series). Iterates with order 3 (such that $g(g(g(x))) = f(x)$, where $f(x)$ is the original function) can be approximated, albeit with decreasing accuracy as the order of the iterate increases.

Using a Fourier series with 32 terms, we obtain accurate results (errors within $10^{-1}$) for fractional iterates of $\sin(x)$ up to order 12.



\begin{figure}[h]
    \centering
    \includegraphics[height=0.4\linewidth]{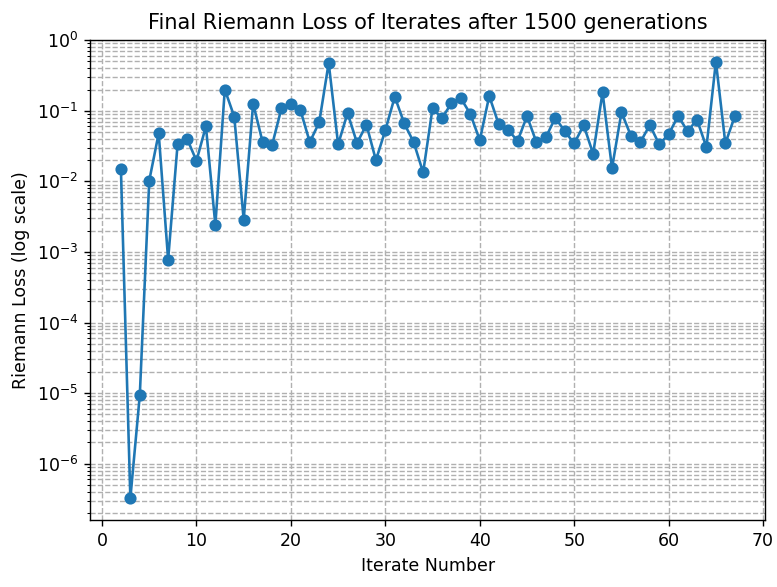}
    \caption{Riemann Losses of higher order genetic approximations for $\sin(x)$ from order 2 to 64 over the interval (-10,10).}
    \label{riemann_loss}
\end{figure}

We see in figure \ref{riemann_loss} that the genetic algorithm technique remains relatively accurate for iterates with a low order but has difficulty producing meaningful results beyond order $1/18$ (the 18th fractional iterate of $\sin(x)$).

\begin{figure}[h]
    \centering
    \includegraphics[height=0.4\linewidth]{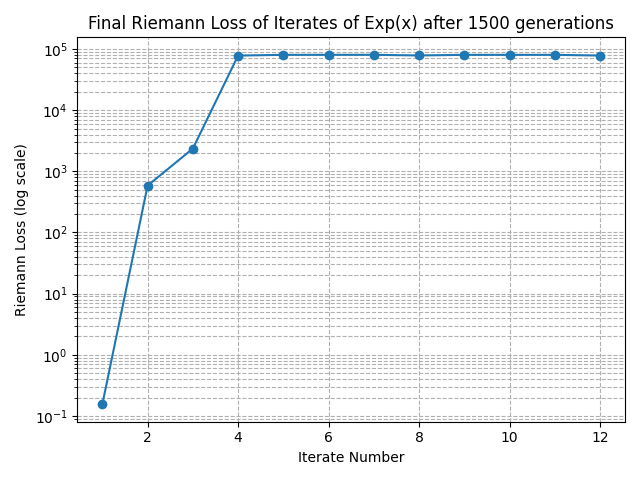}
    \caption{Riemann Losses of higher order genetic approximations for $e^x$ from order 1 to 12 over the interval (-6,6).}
    \label{riemann_loss_exp}
\end{figure}

Another test was run for $e^x$, where power series were evolved using the genetic algorithm to approximate fractional iterates. In figure \ref{riemann_loss_exp}, we see that the error similarly starts low, however, increases for small fractional order iterates due to a lack of terms in the Taylor series and compute power.

For a Taylor series with 32 terms, we obtain accurate results (errors within $10^{4}$) only for fractional iterates of $\sin(x)$ up to order 3, a lower accuracy than that of the $\sin(x)$ genetic algorithm.

We anticipate that with increased computational resources this technique could be used to generate useful fractional iterates up to a meaningfully low fractional order.

\medskip

It is also still possible to obtain iterates of higher order by composing lower order iterates. This technique is outlined as follows:

\begin{enumerate}
    \item Start with a function $f(x)$ and two fractional iterates, $f_\frac{1}{a}(x)$ and $f_\frac{1}{a}(x)$, such that $a$ and $b$ are co-prime, and where $f_n(x)$ denotes a function that is $f(x)$ composed with itself $n$ times.
    \item Compose $f_\frac{1}{a}x)$ with $f_\frac{1}{b}(x)$ to obtain $f_\frac{a+b}{ab}(x)$.
\end{enumerate}

Given that $a$ and $b$ are co-prime, this will result in a new function that has not been generated by any previous technique, as they only generate \emph{harmonic iterates}, meaning those of the form $f_\frac{1}{n}(x)$.

An example is $f_\frac{1}{6}(x)$, which can be obtained by using this technique using $f_\frac{1}{2}(x)$ and $f_\frac{1}{3}(x)$, and then composing the inverse of the result with $f(x)$.

\section*{Conclusion}

In this work, we have explored the challenging functional equation \( f(f(x)) = e^x \) from multiple perspectives. Initially, we established the nonexistence of an elementary closed-form solution by analyzing the derivative \( f'(x) \) and invoking classical results from elementary function theory. This negative result motivated a deeper investigation into the structural properties that any continuous solution must satisfy, including injectivity, monotonicity, and the non-surjectivity of the image of \( f \).

Building on these foundational insights, we presented a known piecewise-defined solution and verified its consistency with the functional equation. Furthermore, by assuming specific functional forms, we derived approximate solutions and connected the problem to the theory of tetration and the super-logarithm, highlighting their relevance to continuous iteration and fractional iterates of the exponential function.

Numerical methods were developed and analyzed, demonstrating that while elementary closed-form solutions do not exist, half iterates and fractional iterates can be approximated with high accuracy. The implementation of iterative chain approximations, Fourier and Taylor series expansions, as well as genetic algorithms, provide practical tools for further exploration.

Future work may focus on refining numerical stability, extending analytic definitions, and exploring broader classes of functional equations with similar iterative structures.

\medskip

Overall, the study illustrates both the depth of theoretical obstacles and the promise of modern computational techniques in advancing our understanding of functional iteration and its applications.

\bibliographystyle{abbrv}

\bibliography{references}

\end{document}